\documentclass[11pt]{article}

\usepackage{amscd,amsmath, amssymb, fancyhdr, epsfig}

\numberwithin{equation}{section}

\newcommand{\version}{version 2.0,\ \   March 2009}

\def\eqref#1{(\ref{#1})}
\newcommand{\goth}{\mathfrak}

\newcommand{\arrow}{{\:\longrightarrow\:}}
\newcommand{\Z}{{\Bbb Z}}
\newcommand{\C}{{\Bbb C}}
\newcommand{\R}{{\Bbb R}}
\newcommand{\Q}{{\Bbb Q}}

\newcommand{\6}{\partial}
\def\1{\sqrt{-1}\:}
\newcommand{\restrict}[1]{{\left|_{{\phantom{|}\!\!}_{#1}}\right.}}
\newcommand{\cntrct}                
{\hspace{2pt}\raisebox{1pt}{\text{$\lrcorner$}}\hspace{2pt}}

\newcommand{\calo}{{\cal O}}

\renewcommand{\tilde}{\widetilde}
\renewcommand{\bar}{\overline}
\renewcommand{\phi}{\varphi}
\renewcommand{\epsilon}{\varepsilon}
\renewcommand{\geq}{\geqslant}
\renewcommand{\leq}{\leqslant}

\newcommand{\Id}{\operatorname{Id}}

\newcommand{\Vol}{\operatorname{Vol}}

\newcommand{\Sym}{\operatorname{Sym}}

\newcommand{\rk}{\operatorname{rk}}

\newcounter{Mycounter}[section]
\newcounter{lemma}[section]
\renewcommand{\thelemma}{{Lemma \thesection.\arabic{lemma}}}
\newcommand{\lemma}{%
    \setcounter{lemma}{\value{Mycounter}}
    \refstepcounter{lemma}
    \stepcounter{Mycounter}
    {\noindent \bf \thelemma:\ }}

\newcounter{claim}[section]
\renewcommand{\theclaim}{{Claim \thesection.\arabic{claim}}}
\newcommand{\claim}{%
    \setcounter{claim}{\value{Mycounter}}
    \refstepcounter{claim}
    \stepcounter{Mycounter}
    {\noindent \bf \theclaim:\ }}

\newcounter{sublemma}[section]
\newcounter{corollary}[section]
\renewcommand{\thecorollary}{{Corollary \thesection.\arabic{corollary}}}
\newcommand{\corollary}{%
    \setcounter{corollary}{\value{Mycounter}}
    \refstepcounter{corollary}
    \stepcounter{Mycounter}
    {\noindent \bf \thecorollary:\ }}

\newcounter{theorem}[section]
\renewcommand{\thetheorem}{{Theorem \thesection.\arabic{theorem}}}
\newcommand{\theorem}{%
    \setcounter{theorem}{\value{Mycounter}}
    \refstepcounter{theorem}
    \stepcounter{Mycounter}
    {\noindent \bf \thetheorem:\ }}

\newcounter{conjecture}[section]
\renewcommand{\theconjecture}{{Conjecture \thesection.\arabic{conjecture}}}
\newcommand{\conjecture}{%
    \setcounter{conjecture}{\value{Mycounter}}
    \refstepcounter{conjecture}
    \stepcounter{Mycounter}
    {\noindent \bf \theconjecture:\ }}

\newcounter{proposition}[section]
\renewcommand{\theproposition}
      {{Proposition \thesection.\arabic{proposition}}}
\newcommand{\proposition}{%
    \setcounter{proposition}{\value{Mycounter}}
    \refstepcounter{proposition}
    \stepcounter{Mycounter}
    {\noindent \bf \theproposition:\ }}

\newcounter{definition}[section]
\renewcommand{\thedefinition}
      {{Definition~\thesection.\arabic{definition}}}
\newcommand{\definition}{%
    \setcounter{definition}{\value{Mycounter}}
    \refstepcounter{definition}
    \stepcounter{Mycounter}
    {\noindent \bf \thedefinition:\ }}

\newcounter{example}[section]
\newcounter{remark}[section]
\renewcommand{\theremark}{{Remark \thesection.\arabic{remark}}}
\newcommand{\remark}{%
    \setcounter{remark}{\value{Mycounter}}
    \refstepcounter{remark}
    \stepcounter{Mycounter}
    {\noindent \bf \theremark:\ }}

\newcounter{problem}[section]
\newcounter{question}[section]
\makeatletter

\renewcommand{\leftmark}{{\scriptsize Hyperk\"ahler SYZ conjecture and semipositive line bundles}}

\@addtoreset{equation}{section} \@addtoreset{footnote}{section}
\makeatother

\def\blacksquare{\hbox{\vrule width 5pt height 5pt depth 0pt}}
\def\endproof{\blacksquare}

\begin{document}
\begin{center}
{\LARGE\bf
Hyperk\"ahler SYZ conjecture and semipositive line bundles\\[4mm]
}

 Misha
Verbitsky

\end{center}

{\small \hspace{0.15\linewidth}
\begin{minipage}[t]{0.7\linewidth}
{\bf Abstract} \\
Let $M$ be a compact, holomorphic symplectic K\"ahler
manifold, and $L$ a non-trivial line bundle admitting a
metric of semi-positive curvature. We show that some power of
$L$ is effective. This result is related to the 
hyperk\"ahler SYZ conjecture, which states that 
such a manifold admits a holomorphic Lagrangian 
fibration, if $L$ is not big.
\end{minipage}
}

\tableofcontents


\section{Introduction}


\subsection{Lagrangian fibrations on hyperk\"ahler manifolds}

Let $M$ be a compact, K\"ahler, holomorphically 
symplectic manifold. By Calabi-Yau theorem,
such a manifold admits a hyperk\"ahler
metric, which is Ricci-flat and necessarily unique
in its K\"ahler class. Throughout
this paper, we shall use 
``hyperk\"ahler'' as 
``compact, K\"ahler, holomorphically 
symplectic''.

Using Bochner vanishing and Berger's
classification of irreducible holonomy,
one proves Bogomolov's decomposition theorem
(\cite{_Besse:Einst_Manifo_}, \cite{_Bogomolov:decompo_}),
stating that any compact hyperk\"ahler
manifold has a finite covering 
$\tilde M \cong T \times M_1 \times ... \times M_n$,
where $T$ is a hyperk\"ahler torus, and 
all $M_i$ are holomorphically symplectic manifolds with 
\begin{equation}\label{_simple_hk_Equation_}
H^1(M_i)=0, \ \  H^{2,0}(M_i)\cong \C
\end{equation}
By Cheeger-Gromoll theorem, a compact Ricci-flat manifold with
$H^1(M_i)=0$ has finite fundamental group. An easy
argument involving the Bochner vanishing would
immediately imply that a hyperkaehler manifold with 
finite fundamental group is in fact simply connected.

A simply connected hyperk\"ahler manifold satisfying 
$H^{2,0}(M)=\C$ is called {\bf simple}. 

In \cite{_Matsushita:fibred_}, D. Matsushita proved the 
following theorem (see also \cite{_Huybrechts:lec_}).

\hfill

\theorem\label{_Matsushita_fibr_Theorem_}
(\cite{_Matsushita:fibred_})
Let $M$ be a simple hyperk\"ahler manifold, and $\pi:\; M \arrow X$
a surjective holomorphic map from $M$ onto a K\"ahler variety $X$.
Assume that $0< \dim X <\dim M$. Then $\dim X = \frac 1 2 \dim M$,
and the fibers of $\pi$ are Lagrangian subvarieties of $M$.

\endproof

Such a map is called {\bf a holomorphic Lagrangian fibration}.
A (real) Lagrangian subvariety $S$ of an $n$-dimensional Calabi-Yau manifold
is called {\bf special Lagrangian} if a holomorphic $(n,0)$-form,
restricted to $S$, is proportional to its Riemannian volume.

From  Calabi-Yau theorem
it follows immediately that
a compact, holomorphically symplectic K\"ahler manifold
admits a triple of complex structures $I, J, K$,
satisfying quaternionic relations (see \cite{_Besse:Einst_Manifo_}).
A complex Lagrangian subvariety of $(M,I)$ is 
special Lagrangian with respect to $J$, which
is clear from the linear algebra. Therefore,
Matsushita's theorem (\ref{_Matsushita_fibr_Theorem_})
gives a way to produce special Lagrangian fibrations
on hyperk\"ahler manifolds.

Holomorphic Lagrangian fibrations are important in Mirror Symmetry.
In \cite{_SYZ:MS_is_T_du_}, Strominger, Yau and Zaslow 
conjectured that Mirror Symmetry of Calabi-Yau manifolds
comes from real Lagrangian fibrations. From arguments
making sense within the framework of string theory,
it occurs that any Calabi-Yau manifold which admits
Mirror Symmetry must apparently admit a 
special Lagrangian fibration, and the dual
fibrations should correspond to the mirror dual
Calabi-Yau manifolds. 

The Strominger-Yau-Zaslow conjecture remains a mystery even now.
For a current survey of SYZ-conjecture, please see \cite{_Gross:SYZ_}.

Examples of special Lagrangian fibrations are 
very rare; indeed, all known examples are derived from
holomorphic Lagrangian fibrations on K3, torus, or other
hyperk\"ahler manifolds.

Existence of holomorphic Lagrangian fibrations on hyperk\"ahler
manifolds is predicted by the Strominger-Yau-Zaslow
interpretation of mirror symmetry. In the weakest form,
the hyperk\"ahler SYZ conjecture is stated as follows.

\hfill

\conjecture\label{_SYZ_weak_Conjecture_}
Let $M$ be a hyperk\"ahler manifold. Then $M$ can be
deformed to a hyperk\"ahler manifold admitting 
a holomorphic Lagrangian fibration.

\hfill

For a more precise form of a 
hyperk\"ahler SYZ conjecture, see Subsection 
\ref{_semiposi_intro_Subsection_}.

The hyperk\"ahler SYZ conjecture is often called
{\bf the Huy\-brechts-\-Sa\-won conjecture}, because it 
was stated in \cite{_Sawon_} and \cite{_Huybrechts:lec_}
(Section 21.4). The same conjecture in 
precise form was stated several years earlier, 
by Hassett and Tschinkel (\cite{_Hassett_Tschinkel:SYZ_conj_}, 
Conjecture 3.8 and Remark 3.12). 

In 1992, in a Harvard lecture, the SYZ conjecture 
for hyperk\"ahler manifolds was stated
by F. A. Bogomolov. However, when I asked 
Bogomolov about the history of this 
conjecture, he said that most likely
it originates in collaboration with A. N. Tyurin
(circa 1985). 

Instead of historically accurate name
``the Tyurin-Bogomolov-Hassett-Tschinkel-Huybrechts-Sawon
conjecure'', we shall call \ref{_SYZ_weak_Conjecture_}
and its precise form \ref{_SYZ_strongest_Conjecture_}
``hyperk\"ahler SYZ conjecture''\footnote{The name
was suggested to the author by Geo Grantcharov.}.

In algebraic geometry, a version of this
conjecture is sometimes called {\em an abundance conjecture},
see \ref{_abundance_Remark_}.

\subsection{Bogomolov-Beauville-Fujiki form on hyperk\"ahler manifolds}
\label{_BBF_Subsection_}

Let $M$ be a simple hyperk\"ahler manifold,
$\dim_\C M = 2n$, $\omega$ its K\"ahler form, and 
$q:\; H^2(M) \times H^2(M)\arrow \R$
a symmetric form on $H^2(M)$ defined by the formula
\begin{equation}\label{_BBF_via_Kahler_Equation_}
   q(\eta_1,\eta_2):=
   \int_X \omega^{2n-2}\wedge \eta_1\wedge\eta_2  
   - \frac{2n-2}{(2n-1)}
   \frac{\int_X \omega^{2n-1}\eta_1 \cdot \int_X\omega^{2n-1}\eta_2}{\int_M\omega^{2n}}
\end{equation}
It is well known (see e. g. \cite{_Beauville_},
\cite{_Verbitsky:cohomo_}, Theorem 6.1,
or \cite{_Huybrechts:lec_}, 23.5, Exercise 30),
that this form is, up to a constant multiplier, independent 
from the choice of complex and K\"ahler structure on $M$,
in its deformation class. The form $q$ is called {\bf
the Bogomolov-Beauville-Fujiki form of $M$}.

Usually, the Bogomolov-Beauville-Fujiki form is
defined as 
\begin{multline*}
   q(\eta,\eta):=
   (n/2)\int_X \eta\wedge\eta  \wedge \Omega^{n-1}
   \wedge \bar \Omega^{n-1} -\\- 
(1-n)\left(\int_X \eta \wedge \Omega^{n-1}\wedge \bar
   \Omega^{n}\right) \left(\int_X \eta \wedge \Omega^{n}\wedge \bar \Omega^{n-1}\right)
\end{multline*}
where $\Omega$ is the holomorphic symplectic form. This definition
(up to a constant multiplier)
is equivalent to the one given above ({\em loc. cit.}).

The Bogomolov-Beauville-Fujiki form largely 
determines the structure of cohomology of $M$,
as implied by the following theorem.

\hfill

\theorem\label{_cohomo_Theorem_}
(see \cite{_Verbitsky:cohomo_},
\cite{_Verbitsky:coho_announce_})
Let $M$ be a simple hyperk\"ahler manifold, $\dim_\C M=2n$,
and $A^{2*}\subset H^{2*}(M)$
the subalgebra in cohomology generated by $H^2(M)$.
Then
\begin{description}
\item[(i)] The natural action of $SO(H^2(M), q)$
on $H^2(M)$ can be extended to a multiplicative action 
on $A^{2*}$. 
\item[(ii)] As an $SO(H^2(M), q)$-representation,
$A^{2i}$ is isomorphic to the symmetric power $\Sym^i(H^2(M))$
for $i \leq n $ and to $\Sym^{2n-i}(H^2(M))$ for $i\geq n$.
\item[(iii)] The properties (i)-(ii) determine the
algebra structure on $A^{2*}$ uniquely.

\item[(iv)] The automorphism group of $A^{2*}$
is isomorphic to $\R^*\times SO(H^2(M), q)$. 
\end{description}

\endproof

\remark
From \ref{_cohomo_Theorem_} (iv),
it follows that the Bogomolov-\-Beauville-\-Fujiki form
is uniquely, up to a constant, determined by topology of $M$.

\hfill

The following theorem was proven by A. Fujiki
in \cite{_Fujiki_}. It follows immediately from the explicit
description of $A^{2*}\subset H^*(M)$ given in \ref{_cohomo_Theorem_}.
This theorem is also sometimes used  as a definition of
Bogomolov-Beauville-Fujiki form.

\hfill

\theorem\label{_Fujiki_formula_Theorem_}
(Fujiki's formula)
Let $M$ be a simple hyperk\"ahler manifold, $\dim_\C M=2n$,
and $a \in H^2(M)$ a non-zero cohomology class. Then
\[
\int_M a^{2n} = q(a,a)^n \lambda,
\]
for some constant $\lambda$ determined by the choice
of Bogomolov-Beauville-Fujiki form.

\endproof

\hfill

From Fujiki's formula and \ref{_cohomo_Theorem_}, Matsushita's theorem 
(\ref{_Matsushita_fibr_Theorem_}) follows quite easily.
Indeed, consider a surjective holomorphic map 
$M \stackrel \pi \arrow X$, with $X$ 
K\"ahler, $0< \dim X < \dim M$,
and let $\omega_X$ be a K\"ahler class
of $X$. Then $(\pi^*\omega_X)^{\dim_\C X}= \pi^*\Vol_X$
is a non-zero positive closed form, hence
its class in $H^*(M)$ is non-zero. Since
$\omega_X^{\dim_\C X+1}$ vanishes identically,
$(\pi^*\omega_X)^{\dim_\C X+1}=0$.
By Fujiki's formula, $q(\pi^*\omega_X,\pi^*\omega_X)=0$.
By \ref{_cohomo_Theorem_}, then,
$(\pi^*\omega_X)^n\neq 0$, where $2n = \dim_\C M$,
and  $(\pi^*\omega_X)^{n+1}=0$. This gives $n =\dim_\C X$.

\subsection{Semipositive line bundles and effectivity}
\label{_semiposi_intro_Subsection_}

\definition
Let $L$ be a line bundle on a compact complex manifold $M$.
Then $L$ is called {\bf semiample} if there exists
a holomorphic map $\pi:\; M \arrow X$
to a projective variety $X$, and $L^N \cong \pi^*(\calo(1))$,
for some $N >0$.

Recall that a cohomology class $\eta \in H^{1,1}(M)$
on a K\"ahler manifold is called {\bf nef} if $\eta$ belongs 
to a closure of the K\"ahler cone of $M$. A line bundle $L$
on $M$ is called {\bf nef} if $c_1(L)$ is nef.

\hfill

From the above argument we obtain that 
holomorphic Lagrangian fibrations
are associated with cohomology classes
$\eta \in H^2(M)$, $q(\eta,\eta)=0$.
Now we can state the hyperk\"ahler SYZ conjecture
in its strongest, most precise form. From now on,
we shall always abbreviate $q(c_1(L), c_1(L))$ 
as $q(L,L)$.

\hfill

\conjecture\label{_SYZ_strongest_Conjecture_}
Let $M$ be simple hyperk\"ahler 
manifold, and $L$ a non-trivial nef bundle
on $M$, with $q(L,L)=0$. Then $L$ is semiample.

\hfill

Notice that the numerical dimension of $L$
is equal to $n=\frac 1 2 \dim_\C M$. By Kawamata's
theorem (\cite{_Kawamata:Pluricanonical_}),
semiampleness of $L$ is implied by
$\kappa(L)=\nu(L)$, where $\kappa$
is Kodaira dimension of $L$, and
$\nu(L)$ is its numerical dimension.
Therefore, to prove \ref{_SYZ_strongest_Conjecture_}
it suffices to show that $H^0(L^x)$
grows as  $x^n$ as $x$ tends to $\infty$.

\hfill

\remark\label{_abundance_Remark_}
When $L$ is the canonical bundle of
a manifold, the equality $\kappa(L)=\nu(L)$ is sometimes called
``the abundance conjecture'' (see
e.g. \cite{_Demailly_Peternell_Schneider:ps-eff_}, 2.7.2). This assertion
is equivalent to the canonical bundle being semiample
(\cite{_Kawamata:Pluricanonical_}). 

\hfill

In \cite{_Matsushita:nef_},
D. Matsushita much advanced this argument.
From his results it follows that $\kappa(L)=\nu(L)$ 
holds if the union of all closed 1-dimensional 
subvarieties $C\subset M$
such that $L\restrict C=\calo_C$ is Zariski
dense in $M$.

If \ref{_SYZ_strongest_Conjecture_}
is true, any nef bundle with
$q(L,L)=0$ should admit a smooth
metric with semi-positive curvature.
The implications of semi-positivity
of $L$ seem to be of independent interest.

Recall that a holomorphic line bundle
is called {\bf effective} if it admits
a non-trivial holomorphic section, and {\bf $\Q$-effective},
if its positive tensor power admits
a section. The main result of the present paper is
the following theorem.

\hfill

\theorem\label{_effe_main_Theorem_}
Let $M$ be a  simple hyperk\"ahler 
manifold, and $L$ a non-trivial nef bundle
on $M$, with $q(L,L)=0$. Assume that $L$ admits a smooth
metric with semi-positive curvature.
Then $L$ is $\Q$-effective.

\hfill

\ref{_effe_main_Theorem_} is proven as follows. 
Using a version of Kodaira-Nakano argument, we 
construct an embedding of $H^i(L)$ to a space of
$L$-valued holomorphic differential forms.
This result is known in the literature
as {\em hard Lefschetz theorem with coefficients
in a bundle} (Section \ref{_hard_Lef_Section_}). 
The holomorphic
Euler characteristic of a line bundle $L$,
denoted as $\chi(L)$, is a polynomial on $q(L,L)$
with coefficients which depend only on the
Chern classes of $M$, as shown by Fujiki
(\cite{_Fujiki_}, 4.12). Therefore,
$\chi(L)= \chi(\calo_M)=n+1$
(\cite{_Besse:Einst_Manifo_}).
This implies existence of
non-trivial $L^N$-valued holomorphic differential forms
on $M$, for any $N>0$ (\ref{_dim_V_k_cohomo_Corollary_}).

To prove \ref{_effe_main_Theorem_},
it remains to deduce that $H^0(L^N)\neq 0$
from existence of $L$-valued holomorphic differential forms
(\ref{_holo_L_valued_L_effe_Theorem_}).
This is done in Section \ref{_L_valued_forms_iplies_effe_Section_}.
We use results of Huybrechts and Boucksom on
duality of pseudo-effective cone and the
modified nef cone on a hyperk\"ahler manifold.
Recall that cohomology class $\eta \in H^{1,1}(M)$
is called {\bf pseudo-effective} if
it can be represented by a closed, positive
current. The pseudo-effective classes form a closed cone
in $H^{1,1}(M)$. A {\bf modified nef} cone is a closure of
the union of all classes $\eta \in H^{1,1}(M)$,
such that for some birational morphism
$\tilde M \arrow M$, $\phi^*\eta$ is nef. In the literature,
the modified nef cone is often called
{\em the movable cone} (this terminology
was introduced by Kawamata) and its interior
the {\em birational K\"ahler cone.}
In \cite{_Huybrechts:cone_} and \cite{_Boucksom_},
Huybrechts and Boucksom prove that 
a dual cone (under the Bogomolov-Beauville-Fujiki
pairing) to the pseudo-effective cone is the
is the modified nef cone of $M$.

Using this duality, and stability of
the tangent vector bundle, we prove that
for any coherent subsheaf $F \subset {\goth T}$,
the class $-c_1(F)$ is pseudo-effective,
where ${\goth T}$ is some tensor power of
$TM$ (\ref{_subsheaves_tange_pseudoeffe_Theorem_}).

This result, together with Boucksom's divisorial
Zariski decomposition (\cite{_Boucksom_}), 
is used to show that any $L$-valued holomorphic 
differential form on $M$ is non-zero in codimension
2, unless $L$ is $\Q$-effective
(\ref{_holo_L_valued_L_effe_Proposition_}).
Then (unless $L^N$ is effective) the above
construction produces infinitely many
sections $s_i \in L^N\otimes \Omega^pM$,
all non-vanishing in codimension 2. Taking the
determinant bundle $D$ of a sheaf generated by
all $s_i$, and using the codimension-2
non-vanishing of $s_i$, we obtain that 
$D \cong L^N$, for some $N> 0$.
By construction, $D$ has non-zero
holomorphic sections. This proves
effectivity of $L^N$, for some $N>0$.

The SYZ-type problem was treated by Campana, Oguiso 
and Peternell in \cite{_COP:non-alge_}, who proved that
a hyperkaehler manifold of complex dimension 4
is either algebraic, has no meromorphic
functions, or admits a holomorphic
Lagrangian fibration. Using an 
argument based on hard Lefschetz theorem with 
coefficients in a bundle, they also proved
the following result. Let $M$ be a hyperkaehler manifold 
of complex dimension $\geq 4$ admitting 
a nef bundle $L$ with $q(L,L)=0$. Then  $M$ admits 
complex subvarieties of dimension at least 2.
This result can be deduced
from \ref{_coho_nonzero_then_effe_Theorem_}
(ii), because a Lelong set of a singular
positive metric on $L$ must be coisotropic,
hence its dimension is $\geq 1/2\dim_\C M$.


\section{Cohomology of semipositive line bundles
and Hard Lefschetz theorem}
\label{_hard_Lef_Section_}


Throughout this section, we shall consider
smooth metrics on line bundles with semipositive
or seminegative curvature. We give simple proofs of
several results which are well known (in different form)
as "Hard Lefschetz theorem with coefficients in a
bundle". This theorem was rediscovered
several times during the 1990-ies: see
\cite{_Enoki:semipositive_}, \cite{_Takegoshi_} and
\cite{_Mourougane_}. We refer to 
\cite{_Demailly_Peternell_Schneider:tan_nef_},
where the multiplier ideal version of this result
is stated and proven.

\subsection{Harmonic forms with coefficients in a
  semipositive line bundle}

\definition
Let $M$ be a complex manifold, and 
$L$ a holomorphic Hermitian line bundle.
We say that $L$ is {\bf semipositive (seminegative)}
if its curvature is a positive or negative
(but not necessarily positive definite)
$(1,1)$-form. 
Let $K$ denote the canonical bundle of $M$.
The standard proof of Kodaira-Nakano theorem
can be used to show that for any positive bundle $L$, one
has $H^i(L \otimes K)=0$ for all $i>0$,
can be generalized for semipositive bundles.
In semipositive case, we obtain that any non-zero cohomology
class $\eta \in H^i(L \otimes K)=0$ corresponds to 
a non-zero holomorphic $(K \otimes L)$-valued 
$i$-form on $M$. 

\hfill

Let $B$ be a Hermitian line bundle on a 
K\"ahler manifold, and
\[ \bar \6:\; \Lambda^{p,q}(M)\otimes B \arrow \Lambda^{p,q+1}(M)\otimes B,\\
\6:\; \Lambda^{p,q}(M)\otimes B \arrow \Lambda^{p+1,q}(M)\otimes B
\]
 the $(0,1)$ and $(1,0)$-parts of the Chern connection,
$\6^*$ and $\bar\6^*$ the Hermitian adjoint operators,
and $\Delta_{\bar\6}:= \bar\6\bar\6^*+\bar\6^*\bar\6$.
Forms which lie in the kernel of this operator are
called {\bf $\Delta_{\bar\6}$-harmonic}.
From the usual arguments one obtains that
the space of $\Delta_{\bar\6}$-harmonic 
$(0,k)$-forms is identified with 
the holomorphic cohomology $H^k(B)$.

\hfill

\lemma\label{_harmo_semipo_Lemma_}
Let $(M, I, \omega)$ be a compact K\"ahler manifold,
and $B$ a seminegative holomorphic line bundle. 
Consider a $\Delta_{\bar\6}$-harmonic 
$(0,k)$-form $\eta \in \Lambda^{0,k}(M)\otimes B$.
Then $\6\eta=0$. Moreover, $\eta\wedge\Theta=0$,
where $\Theta\in \Lambda^{1,1}(M)$ is the curvature of $B$.

\hfill

{\bf Proof:} Let 
\[ L:\; \Lambda^{p,q}(M)\otimes B\arrow \Lambda^{p+1,q+1}(M)\otimes B\]
be the Hodge operator of multiplication by the K\"ahler form
$\omega$, 
\[ \Lambda:\; \Lambda^{p,q}(M)\otimes B\arrow \Lambda^{p-1,q-1}(M)\otimes B\]
its Hermitian adjoint. The Kodaira
identities are well-known (\cite{_Griffi_Harri_}),
\[
[L, \6^*] = \1 \bar\6,\ \ [\Lambda, \6] = -\1 \bar\6^* 
\]
From these identities one obtains (as usual)
\[
\Delta_{\bar\6}- \Delta_{\6} = - [L_\Theta, \Lambda],
\]
where $L_\Theta$ is an operator of a multiplication by $\Theta$.
Choose an orthonormal frame $\xi_1, ... \xi_n, \bar \xi_1,
..., \bar \xi_n$, such that 
$\Theta = - \sum \alpha_i \xi_i\wedge \bar\xi_i$,
and $\alpha_i$ are non-negative real numbers.
For any $(0,k)$-form $e= \bar\xi_{i_1} \wedge ... \wedge \bar\xi_{i_k}$,
we have
\begin{equation}\label{_commu_on_monom_Equation_}
[L_\Theta, \Lambda]e=\sum \alpha_{j_p} e,
\end{equation}
where the sum is taken over all $\alpha_{j_p}$
with $j_p\notin\{i_1, i_2, ..., i_k\}$. Therefore,
\[
\Delta_{\bar\6} = \Delta_{\6} + A,
\]
where the operator $A(e)= \sum \alpha_{j_p} e$
is positive and self-adjoint on $\Lambda^{0,k}(M)\otimes B$.
This gives, for any $\eta \in \ker \Delta_{\bar\6}$, that
$\Delta_{\6}(\eta)=0$ and $A(\eta)=0$. From
\eqref{_commu_on_monom_Equation_} it is clear 
that $A(\eta)=0$ if and only if $\eta$ is a sum of
monomials $e= \bar\xi_{i_1} \wedge ... \wedge\bar\xi_{i_k}$
containing all $\bar\xi_{i_p}$ with 
$\alpha_{i_p}\neq 0$. This is equivalent to 
$e \wedge \Theta=0$.

We obtained that 
\[ \Delta_{\bar\6}\eta=0 \Leftrightarrow \Delta_{\6}\eta=0 
\text{\ \ and\ \ } \eta \wedge \Theta=0.
\]
However, 
\[ 
  (\Delta_{\6}\eta, \eta) = (\6\eta, \6\eta) + (\6^*\eta, \6^*\eta),
\]
hence $\Delta_{\6}\eta=0$ implies $\6\eta=0$.\footnote{The
  converse is also easy to see; see
  \ref{_harmo_via_holo_Proposition_}.} 
We proved 
\ref{_harmo_semipo_Lemma_}. \endproof

\hfill

Let $B$ be a holomorphic Hermitian line bundle,
and $\nabla$ its Chern connection. Denote by $(\bar B, \bar \nabla)$
 the same bundle with the opposite complex structure
and the same connection. The $(0,1)$-part of $\bar \nabla$
is complex conjugate to $\nabla^{1,0}$, hence
$\left(\bar \nabla^{0,1}\right)^2=0$. Therefore, as follows
from the vector bundle version of Newlander-Nirenberg 
theorem, the operator $\bar \nabla^{0,1}$ defines
a holomorphic structure on $\bar B$. 
This allows one to consider $\bar B$ as a holomorphic
vector bundle. 

\hfill

\claim 
In these assumptions, $\bar B$ is isomorphic 
to $B^*$, as a holomorphic Hermitian vector bundle.

\hfill

{\bf Proof:} The Hermitian metric 
$h:\; B \times \bar B\arrow \C$
is non-degenerate, complex linear, and
preserved by the connection, hence
the corresponding pairing identifies
$\bar B$ as a bundle with connection with $B^*$.
\endproof

\hfill

\proposition\label{_harmo_via_holo_Proposition_}
 Let $(M, I, \omega)$ be a compact K\"ahler manifold,
and $B$ a seminegative holomorphic line bundle.
Consider a form $\eta \in \Lambda^{0,k}(M)\otimes B$,
and let $\bar\eta \in \Lambda^{k,0}(M)\otimes \bar B$
be its complex conjugate. Then the following
statements are equivalent.

\begin{description}
\item[(i)] $\Delta_{\6}\eta=0$
\item[(ii)] $\bar\eta$ is a holomorphic section
of 
\[  
   \Lambda^{k,0}(M)\otimes \bar B\cong \Omega^k M \otimes B^*,
\]
and, moreover, $\bar\eta\wedge \Theta=0$, where $\Theta$
is a curvature form of $B$. 
\end{description}

{\bf Proof:} \ref{_harmo_via_holo_Proposition_}
follows immediately from the same argument
as used to prove \ref{_harmo_semipo_Lemma_}.
We have
\[
\Delta_{\bar\6} = \Delta_{\6} + A
\]
\eqref{_commu_on_monom_Equation_}, and $A$ is positive self-adjoint,
hence $\Delta_{\bar\6}\eta=0$ is equivalent to
\[ \6\eta=\6^*\eta= \eta\wedge \Theta=0.\]
This is obviously equivalent to
\[
\bar\6 \bar\eta = \bar\6^* \bar\eta=\bar\eta\wedge \Theta=0.
\]
However, $\bar\eta$ is a $\bar B$-valued $(k,0)$-form,
hence $\bar\6 \bar\eta$ means that it is holomorphic.
Then $\bar\6^* \bar\eta$ is automatic. We obtained that
$\Delta_{\bar\6}\eta=0$ is equivalent to 
$\bar\6 \bar\eta = \bar\eta\wedge \Theta=0$.
\ref{_harmo_via_holo_Proposition_} is proven.
\endproof

\subsection{Cohomology vanishing for semipositive line bundles}

\theorem\label{_coho_and_holo_Serre_Theorem_}
Let $(M, I, \omega)$ be a compact K\"ahler manifold, $\dim_\C M=n$,
$K$ its canonical bundle, and $B$ a Hermitian holomorphic
line bundle on $M$. Assume that $B^*\otimes K$ is
seminegative, and denote its curvature by $\Theta$. 
Then the following spaces are naturally isomorphic,
for all $k$. 
\begin{description}
\item[(i)] The space $V_k$ of holomorphic forms
$\eta \in \Lambda^{k,0}(M)\otimes B^* \otimes K$ satisfying
$\eta\wedge \Theta=0$.
\item[(ii)] $H^{n-k}(B)^*$.
\end{description}

{\bf Proof:} 
By \ref{_harmo_via_holo_Proposition_},
the space $V$ is isomorphic to $H^k(B^* \otimes K)$.
By Serre's duality, $H^k(B^* \otimes K)$ is dual to $H^{n-k}(B)$.
\endproof

\hfill

\corollary\label{_coho_and_har_CY_Corollary_}
Let $(M, I, \omega)$ be a compact K\"ahler manifold
with trivial canonical bundle, $\dim_\C M=n$, $B$ 
a semipositive line bundle on $M$, and $\Theta$
its curvature. Then the following vector spaces are
naturally isomorphic.
\begin{description}
\item[(i)] The space $V_k$ of holomorphic forms
$\eta \in \Lambda^{k,0}(M)\otimes B^* $ satisfying
$\eta\wedge \Theta=0$.
\item[(ii)] $H^{n-k}(B)^*$.
\end{description}

\endproof

\hfill

Using these arguments for $L^2$-cohomology as in 
\cite{_Demailly:L^2_}, Chapter 5, we could obtain a Nadel vanishing
version of \ref{_coho_and_har_CY_Corollary_}. 
To avoid stating the necessary results and
definitions, we refer directly to \cite{_Demailly_Peternell_Schneider:tan_nef_}.

\hfill

\theorem\label{_Nadel_holo_coho_CY_Theorem_}
(\cite{_Demailly_Peternell_Schneider:tan_nef_}, Theorem 2.1.1)
Let $(M, I, \omega)$ be a compact K\"ahler manifold, $\dim_\C M=n$, 
$K$ its canonical bundle,
and $L$ a holomorphic line bundle on $M$ equipped with a singular
Hermitian metric $h$. Assume that the 
curvature $\Theta$ of $L$ is a positive current on $M$,
and denote by ${\cal I}(h)$ the corresponding multiplier ideal
(Section \ref{_Multiplie_Section_}).
Then the wedge multiplication operator $\eta \arrow \omega^i \wedge \eta$
induces a surjective map
\[
H^0(\Omega^{n-i}M\otimes  L\otimes {\cal I}(h))\stackrel {\omega^i \wedge
\cdot}\arrow H^i(K  \otimes  L\otimes {\cal I}(h)).
\]
\endproof

\hfill

\ref{_coho_and_har_CY_Corollary_} immediately
leads to some interesting results about
semipositive bundles on hyperk\"ahler manifolds. 

\hfill

\corollary\label{_dim_V_k_cohomo_Corollary_}
Let $(M,I)$ be a simple hyperk\"ahler
manifold, $\dim_\C M =2n$, and $L$ a non-trivial nef bundle which satisfies
$q(L,L)=0$. Assume that $L$ admits a Hermitian metric with
semipositive curvature form $\Theta$. Then

\begin{description}
\item[(i)] $H^i(L)=0$, for all $i>n$.
\item[(ii)] Denote by $V_k$ the space $H^0(\Omega^k\otimes  L)$ 
of holomorphic $L$-valued forms $\eta \in
\Lambda^{k,0}(M)\otimes L$ satisfying $\eta\wedge\Theta=0$. 
Then $\dim V_k = \dim H^{2n-k}(L)$.
\item[(iii)] $\sum (-1)^k \dim V_k = n+1$.
\end{description}

{\bf Proof:} \ref{_dim_V_k_cohomo_Corollary_} (i) follows
from \cite{_Verbitsky:qD_}, Theorem 1.7, and Proposition~6.4. 
\ref{_dim_V_k_cohomo_Corollary_} (ii) is a restatement of
\ref{_coho_and_har_CY_Corollary_}.  
\ref{_dim_V_k_cohomo_Corollary_} (ii) follows
from \cite{_Fujiki_}, 4.12. Indeed,
as Fujiki has shown, the holomorphic Euler characteristic 
$\chi(L)$ of $L$ is expressed through $q(L,L)$. Therefore,
it is equal to the holomorphic
Euler characteristic $\chi(\calo_M)$ of the trivial sheaf $\calo_M$.
However, 
$\chi(\calo_M) = n+1$, as follows from Bochner's vanishing
theorem (see e.g. \cite{_Besse:Einst_Manifo_}).
\endproof

\hfill

\remark
In assumptions of \ref{_dim_V_k_cohomo_Corollary_},
we established existence of at least $n+1$ linearly
independent holomorphic $L$-valued forms on $M$.

\hfill

\remark
In \cite{_Demailly_Peternell_Schneider:ps-eff_},
\ref{_Nadel_holo_coho_CY_Theorem_} was used to
obtain a weak form of Abundance Conjecture for manifolds $M$
with pseudo-effective canonical bundle $K_M$ 
admitting a singular metric with algebraic singularities. 
It was shown here (Theorem 2.7.3) that either such a manifold
admits a non-trivial holomorphic differential form,
or $H^0(\Omega^*(M)\otimes K_M^{\otimes m}$
is non-zero for infinitely many $m>0$.


\section{Stability and $L$-valued holomorphic forms}
\label{_L_valued_forms_iplies_effe_Section_}


The main result of this section is the following
theorem, proven in Subsection \ref{_holo_forms_then_effe_Subsection_}.

\hfill

\theorem\label{_holo_L_valued_L_effe_Theorem_}
Let $M$ be a compact hyperk\"ahler manifold, and
$L$ a nef line bundle satisfying $q(L, L)=0$.
Assume that $M$  admits a non-trivial $L^k$-valued
holomorphic differential form, for infinite number
of $k\in \Z^{>0}$.\footnote{This would follow
immediately from $\chi(L^k)=n+1$, if $L$ admits a 
semi-positive metric, by the result of Fujiki
which ismentioned in \ref{_dim_V_k_cohomo_Corollary_}.}
  Then $L$ is $\Q$-effective.

\subsection{Stability and Yang-Mills connections}

We remind some standard facts about the
Kobayashi-Hitchin correspondence, proven
by Donaldson and Uhlenbeck-Yau 
(\cite{_Uhle_Yau_}, \cite{_Lubcke_Teleman_})

\hfill

\definition\label{_degree,slope_destabilising_Definition_} 
Let $F$ be a coherent sheaf over
an $n$-dimensional compact K\"ahler manifold $M$. We define
{\bf the degree} $\deg(F)$ as
\[ 
   \deg(F)=\int_M\frac{ c_1(F)\wedge\omega^{n-1}}{vol(M)}
\] 
and $\text{slope}(F)$ as
\[ 
   \text{slope}(F)=\frac{1}{\text{rank}(F)}\cdot \deg(F). 
\]

Let $F$ be a torsion-free coherent sheaf on $M$
and $F'\subset F$ a proper subsheaf. Then $F'$ is 
called {\bf destabilizing subsheaf} 
if $\text{slope}(F') \geq \text{slope}(F)$

A coherent sheaf $F$ is called {\bf 
 stable}
\footnote{In the sense of Mumford-Takemoto}
if it has no destabilizing subsheaves. 
A coherent sheaf $F$ is called {\bf 
polystable} if it is a direct sum of stable sheaves of the same slope.
 
\hfill

\definition \label{_Yang-Mills_Definition_}
Let $B$ be a holomorphic Hermitian   bundle over a
K\"ahler manifold $M$, $\nabla$ its Chern connection, and 
$\Theta\in\Lambda^{1,1}\otimes End(B)$ its curvature.
The Hermitian metric on $B$ and the connection $\nabla$
defined by this metric are called {\bf 
Yang-Mills} if 
\[
   \Lambda(\Theta)=constant\cdot \Id\restrict{B},
\]
where $\Lambda$ is a Hodge operator and $\Id\restrict{B}$ is 
the identity endomorphism which is a section of $End(B)$.

\hfill

\theorem \label{_UY_Theorem_} 
(Uhlenbeck-Yau)
Let B be an indecomposable
holomorphic bundle over a compact K\"ahler manifold. Then $B$ admits
a Hermitian 
Yang-Mills connection if and only if it is 
stable. Moreover, the Yang-Mills 
connection is unique, if it exists.
 
{\bf Proof:} \cite{_Uhle_Yau_}. \endproof

\hfill

\remark
Any tensor power of a Yang-Mills bundle
is again Yang-Mills. This implies that
a tensor power of a polystable bundle is again
polystable. Notice that this result
follows from \ref{_UY_Theorem_}.

\hfill

\remark\label{_Einstein_Yang_Mills_Remark_}
Given a K\"ahler-Einstein manifold (e.g. a 
Calabi-Yau, or a hyperk\"ahler manifold),
its tangent bundle is manifestly Yang-Mills
(the curvature condition $\operatorname{Ric}(M)=const$
is equivalent to the Yang-Mills condition,
as follows from a trivial linear-algebraic
argument; see \cite{_Besse:Einst_Manifo_} for details). 
Therefore, $TM$ is polystable, for any 
Calabi-Yau manifold.

\subsection{The birational K\"ahler cone}

\definition
(\cite{_Huybrechts:cone_}, see also \cite{_Boucksom_})
Let $(M,\omega)$ be a compact K\"ahler manifold, and
$\{(M_\alpha, \phi_\alpha)\}$ the set of all compact
manifolds equipped with a birational
morphism $\phi_\alpha:\; M_\alpha\arrow M$. 
 Let ${\cal MN}(M)\subset H^{1,1}(M)$ be the closure of 
the set of all classes $\eta\in H^{1,1}(M)$ such that for some
$(M_\alpha, \phi_\alpha)$, the pullback 
$\phi_\alpha^*\eta$ is a K\"ahler class on $M_\alpha$.
The set ${\cal MN}(M)$ is called
{\bf the modified nef cone}, and its
interior part {\bf the birational
K\"ahler cone} of $M$, or {\bf the modified K\"ahler cone}.

\hfill

\definition
Let $M$ be a compact K\"ahler manifold, and $\eta\in H^{1,1}_\R(M)$
a real $(1,1)$-class which can be represented by 
a positive, closed (1,1)-current. Then $\eta$
is called {\bf pseudoeffective}.

\hfill

\theorem\label{_MN_dual_to_effe_Theorem_}
Let $M$ be a compact hyperk\"ahler manifold, and 
${\cal MN}(M)\subset H^{1,1}_\R(M)$ its birational 
nef cone. Denote by ${\cal P}(M)\subset H^{1,1}_\R(M)$
its pseudoeffective cone. Then ${\cal P}(M)$ is dual
to ${\cal MN}(M)$ with respect to Bogomolov-Beauville-Fujiki
pairing. 

\hfill

{\bf Proof:} \cite{_Huybrechts:cone_}, Proposition 4.7, 
 \cite{_Boucksom_}, Proposition 4.4. \endproof

\hfill

\remark\label{_MK_Huybrechts_de_Remark_}
Let $M$ be a hyperk\"ahler manifold, and 
$M_\alpha$ another hyperk\"ahler manifold, which is
birationally equivalent to $M$. It is well-known that
there is a natural isomorphism 
$H^{1,1}(M)\cong H^{1,1}(M_\alpha)$.
In \cite{_Huybrechts:cone_}, 
Huybrechts defined the birational K\"ahler cone as a 
inner part of a cone obtained as 
a closure of a union of all K\"ahler cones ${\cal K}(M_\alpha)$,
for all hyperk\"ahler manifolds $M_\alpha$ 
birationally equivalent to $M$. 
This definition is equivalent to the one
given above, as shown in \cite{_Boucksom_}.

\hfill

\theorem\label{_subsheaves_tange_pseudoeffe_Theorem_}
Let $M$ be a compact hyperk\"ahler manifold, ${\goth T}$  
a tensor power of a tangent bundle (such as a bundle
of holomorphic forms), and $E \subset {\goth T}$ a coherent subsheaf
of ${\goth T}$. Then the class $-c_1(E)\in H^{1,1}_\R(M)$
is pseudoeffective.

\hfill

\remark 
In \cite{_Campana_Peternell:Geom_Stab_}, Theorem 0.3,
Campana and Peternell prove that for any projective manifold $X$,
and any surjective map $(\Omega^1X)^{\otimes m}\arrow S$,
with $S$ a torsion-free sheaf, $\det S$ is pseudo-effective,
unless $X$ is uniruled. This general (and beautiful)
result easily implies \ref{_subsheaves_tange_pseudoeffe_Theorem_},
if $M$ is projective. However, its proof is quite difficult,
and does not work for non-algebraic K\"ahler manifolds.

\hfill

{\bf Proof of \ref{_subsheaves_tange_pseudoeffe_Theorem_}:} 
Since $M$ is hyperk\"ahler, $TM$ is 
a Yang-Mills bundle (\ref{_Einstein_Yang_Mills_Remark_}), of slope 0. 
Therefore, its tensor power ${\goth T}$
is also a Yang-Mills bundle. Since
a Yang-Mills bundle is polystable, we have
\begin{equation}\label{_c_1_E_wedge_omega_negati_Equation_}
\int_M c_1(E) \wedge \omega^{\dim_\C M-1}\leq 0,
\end{equation}
for any K\"ahler form $\omega$ on $M$.
Using the formula \eqref{_BBF_via_Kahler_Equation_}, 
we can express the integral
\eqref{_c_1_E_wedge_omega_negati_Equation_} in terms
of the Bogomolov-Beauville-Fujiki form, obtaining
that $\int_M c_1(E) \wedge \omega^{\dim_\C M-1}$
is proportoinal to $q(c_1(E), \omega)$, with positive 
coefficient. Therefore,
\eqref{_c_1_E_wedge_omega_negati_Equation_} 
holds if and only if $-c_1(E)$ lies in the dual
K\"ahler cone ${\cal K}^*(M)$. 

Consider a hyperk\"ahler manifold $M_\alpha$
which is birationally equivalent to $M$, 
let $M_1 \subset M \times M_\alpha$ be the correspondence
defining this birational equivalence, and 
$\pi, \sigma:\; M_1 \arrow M, M_\alpha$ 
the corresponding projection maps. 
Since the canonical class of $M$, $M_\alpha$ is trivial,
the birational equivalence 
$M \stackrel{\phi_\alpha}\arrow M_\alpha$
is an isomorphism outside of codimension 2.
This allows one to identify $H^2(M)$ and $H^2(M_\alpha)$.
Let $Z\subset M_\alpha$ be a set where
$\phi_\alpha$ is not an isomorphism.
Outside of $Z$, the sheaf $\sigma_* \pi^* {\goth T}$ is 
isomorphic to a similar tensor power ${\goth T}_\alpha$
on $M_\alpha$. Therefore, the reflexization
$(\sigma_* \pi^* {\goth T})^{**}$ is naturally 
isomorphic to ${\goth T}_\alpha$. 

Consider the  sheaf $E_\alpha:= (\sigma_* \pi^* E)^{**}$.
Outside of $Z$, $\sigma_* \pi^* E$ is naturally embedded to 
$\sigma_* \pi^* {\goth T}$. Therefore, the corresponding
map of reflexizations is also injective, and we
may consider $E_\alpha$ as a subsheaf of ${\goth T}_\alpha$.
Since $\phi_\alpha$ is an isomorphism outside of codimension 2,
$c_1(E)=c_1(E_\alpha)$. Applying \eqref{_c_1_E_wedge_omega_negati_Equation_}
again, we find that $-c_1(E)$ lies in the dual
K\"ahler cone ${\cal K}^*(M_\alpha)$.
We have shown that
\[
- c_1(E) \in \bigcap_\alpha {\cal K}^*(M_\alpha),
\] 
for all hyperk\"ahler birational modifications 
$M_\alpha$ of $M$. From \ref{_MK_Huybrechts_de_Remark_},
we obtain that $\bigcap_\alpha {\cal K}^*(M_\alpha)$
is the dual cone to the birational K\"ahler cone of $M$.
However, the birational K\"ahler cone is dual to 
the pseudoeffective cone, as follows
from \ref{_MN_dual_to_effe_Theorem_}. We have shown
that $- c_1(E)$ is pseudoeffective.
\ref{_subsheaves_tange_pseudoeffe_Theorem_} is proven.
\endproof

\subsection{Zariski decomposition and $L$-valued holomorphic forms}

The following easy lemma 
directly follows from the Hodge index theorem

\hfill 

\lemma\label{_Hodge_index_Lemma_}
Let $M$ be a compact hyperk\"ahler manifold, 
$\eta\in H^{1,1}(M)$ a nef class satisfying
$q(\eta,\eta)=0$, and $\nu\in H^{1,1}(M)$
a class satisfying $q(\eta, \nu)=0$ and 
$q(\nu, \nu) \geq 0$.
Then $\eta$ is proportional to $\nu$.

\hfill

{\bf Proof:} Suppose $\eta$ is not proportional
to $\nu$. Let $\nu = k\eta + \nu'$,
where $\nu'$ is orthogonal to $\eta$, and $W\subset H^{1,1}_\R(M)$
be a 2-dimensional subspace generated by $\nu', \eta$.
By the Hodge index theorem, the form $q$ on $H^{1,1}_\R(M)$
has signature $(+, -,-,-,-, ...)$ (\cite{_Beauville_}),
hence $q(\nu', \nu') <0$. However, 
$q(\nu, \nu) = q(\nu', \nu')$, because
$q(\eta, \nu')= q(\eta, \nu)=0$. We obtained
a contradiction, proving \ref{_Hodge_index_Lemma_}.
\endproof

\hfill

\remark\label{_Hodge_index_Remark_}
In the sequel, this lemma is applied to the following
situation: $\eta$ is a nef class, satisfying $q(\eta, \eta)=0$, 
and $\nu$ a modified nef class. Then $q(\eta, \nu)=0$
implies that $\eta$ is proportional to $\nu$.

\hfill

\proposition\label{_holo_L_valued_L_effe_Proposition_}
Let $M$ be a compact hyperk\"ahler manifold, 
$L$ a nef line bundle satisfying $q(L, L)=0$,
${\goth T}$ some tensor power of a tangent bundle, 
and  $\gamma \in {\goth T}\otimes L$ a non-zero holomorphic
section.  Consider the zero divisor $D$ 
of $\gamma$ (the sum of all divisorial
components of the zero set of $\gamma$
with appropriate multiplicities). Assume that 
$L$ is not $\Q$-effective. Then $D$ is trivial.

\hfill

{\bf Proof:} Let $L_0$ be a rank 1 subsheaf of ${\goth T}$
generated by $\gamma\otimes L^{-1}$. 
By \ref{_subsheaves_tange_pseudoeffe_Theorem_},
$\nu:= - c_1(L_0)$ is pseudoeffective. Clearly,
$[D]= c_1(L\otimes L_0)$. Therefore, $c_1(L)=[D]+\nu$.
To prove \ref{_holo_L_valued_L_effe_Proposition_}
we are going to show that $\nu$ is proportional
to $c_1(L)$.

Since $c_1(L)$ is a limit of K\"ahler classes,
we have $q(L, D)\geq 0$ and $q(L, \nu)\geq 0$.
Since $0= q(L,L) = q(L, D)+ q(L,\nu)$,
this gives
\[
q(L,D) = q(L, \nu) =0.
\]

In \cite{_Boucksom_}, Proposition 3.10, S. Boucksom
has constructed {\em the Zariski decomposition}
for pseudoeffective classes, showing that any
pseudoeffective class $\nu$ can be decomposed
as $\nu = \nu_0 + \sum \lambda_i [D_i]$,
where $\lambda_i$ are positive numbers,
$D_i$ exceptional divisors, and $\nu_0$ is
a modified nef class.
On a hyperk\"ahler manifold, the numbers $\lambda_i$ are
rational, if $\nu$ is a rational class (\cite{_Boucksom_},
Corollary 4.11).

Since $\eta:=c_1(L)$ is nef, it is obtained as a limit of 
K\"ahler classes, hence $q(\eta, D_i) \geq 0$, 
and $q(\eta, \nu_0)\geq 0$. Therefore,  $q(L, \nu) =0$ 
implies that $q(\eta, \nu_0)=0$  and $q(\eta, D_i)=0$.
By \ref{_Hodge_index_Lemma_}, a modified 
nef class $\nu_0$ which satisfies 
$q(\eta, \nu_0)=0$ is proportional to $\eta$
(see \ref{_Hodge_index_Remark_}). 
Therefore, $\nu = \lambda c_1(L)  + \sum \lambda_i [D_i]$,
where $\lambda\geq 0$.
We obtain
\begin{equation}\label{_L_via_nu_and_Zar_Equation_}
c_1(L) = \lambda c_1(L)  + \sum \lambda_i [D_i]+ [D].
\end{equation}
From \eqref{_L_via_nu_and_Zar_Equation_},
we immediately infer that $(1-\lambda)^{-1}c_1(L)$ is effective,
unless $\lambda=1$, $[D]$ is trivial and all $\lambda_i$ vanish.
\endproof

\hfill

\remark
Using the terminology known from algebraic geometry,
\ref{_Hodge_index_Lemma_} can be rephrased by saying
that a nef class $\eta \in H^{1,1}(M)$ which satisfies 
$q(\eta, \eta)=0$ generates an extremal ray in the
nef cone. Then \ref{_holo_L_valued_L_effe_Proposition_}
would follow from already known arguments
(see also {\em e.~g.} \cite{_Campana_Peternell:Geom_Stab_},
Corollary 1.12).

\subsection{$L$-valued holomorphic forms
on hyperk\"ahler manifolds}
\label{_holo_forms_then_effe_Subsection_}

Now we can prove \ref{_holo_L_valued_L_effe_Theorem_}.
Let $M$ be a compact hyperk\"ahler manifold, and
$L$ a nef line bundle satisfying $q(L, L)=0$.
Assume that $M$  admits a non-trivial $L^k$-valued
holomorphic differential form, for infinite number
of $k\in \Z^{>0}$. We have to show that
$L^{\otimes N}$ is effective, for some $N>0$.

Suppose that $L^{\otimes k}$ is never effective.
Then, by \ref{_holo_L_valued_L_effe_Proposition_},
any non-zero section of ${\goth T}\otimes L^{\otimes k}$
is non-zero outside of codimension 1. 

Let $E=\bigoplus _i \Omega^iM$ be the bundle of all
differential forms, and $E_k\subset E$ its subsheaf
generated by global sections of 
$E\otimes L^{\otimes i}$, $i=1, ..., k$.
Since $E_1 \subset E_2 \subset ...$, this
sequence stabilizes. Let $E_\infty\subset E$
be its limit, $\rk E_\infty=r$. Choose an
$r$-tuple  $\gamma_1 \in E\otimes L^{\otimes i_1}, ..., 
\gamma_r \in E\otimes L^{\otimes i_r}$, 
of linearly independent sections of $E_\infty$. Then
the top exterior product $\gamma_1\wedge ... \wedge \gamma_r$
is a section of $\det E_\infty \otimes L^{\otimes I}$,
where $I = \sum_{k=1}^r i_k$. The rank 1 sheaf
$\det E_\infty \subset \Lambda^r E$
is a subsheaf in a tensor bundle
\[ \Lambda^r E = \Lambda^r\left(\bigoplus _i \Omega^iM\right),
\]
hence any section of $\det E_\infty \otimes L^{\otimes I}$ 
is non-zero in codimension 2 
(\ref{_holo_L_valued_L_effe_Proposition_}).
Therefore, $(\det E_\infty)^*\cong L^{\otimes I}$.

There is infinite number of
$\gamma_k \in E\otimes L^{\otimes i_k}$ to choose,
hence for appropriate choice of 
$\{\gamma_k\in E\otimes L^{\otimes i_k}\}$, the number 
$I = \sum_{k=1}^r i_k$ can be chosen as big as we wish.
Therefore, the isomorphism $(\det E_\infty)^*\cong L^{\otimes I}$
cannot hold for most choices of the set $\{\gamma_k\}$.
We came to contradiction, proving effectivity
of $L^{\otimes k}$ for some $k >0$.
\ref{_holo_L_valued_L_effe_Theorem_} is proven.
\endproof

\section{Multiplier ideal sheaves}
\label{_Multiplie_Section_}

Let $\psi:\; M \arrow [-\infty, \infty[$ 
be a plurisubharmonic function on
a complex $n$-dimensional manifold $M$, and $Z:= \psi^{-1}(-\infty)$.
Recall that such a subset is called 
{\bf a pluripolar set}. It is easy to check
that a complement to a pluripolar set is 
open and dense.

By definition, a {\bf singular metric} on 
a line bundle $L$ is a metric of form 
$h=h_0e^{-2\psi}$, where $\psi$ is a 
locally integrable function, defined 
outside of a closed pluripolar set.

A function is called {\bf quasi-plurisubharmonic}
if it can be locally expressed as a sum of a smooth function
and a plurisubharmonic function.

Let $L$ be a nef bundle on a compact K\"ahler
manifold $M$. Then $c_1(L)$ is a limit of a K\"ahler classes
$\{\omega_i\}$, which are uniformly bounded.
Since the set of positive currents is relatively
compact, the sequence $\{\omega_i\}$ has a limit $\Xi$,
which is a closed, positive current on $M$, representing
$c_1(L)$. Consider a smooth, closed form $\theta$,
representing $c_1(L)$. Using $\6\bar\6$-lemma for
currents, we may assume that $\Xi-\theta=\6\bar\6\psi$,
where $\psi$ is a 0-current, that is, an
$L^1$-integrable distribution. Clearly, $\psi$
is quasi-plurisubharmonic; in particular, 
$\psi$ is upper semi-continuous and locally 
bounded outside of a pluripolar set.

Let $h_0$ be
a Hermitian metric on $L$ such that $\theta$
is its curvature (such a metric always exists
by $\6\bar\6$-lemma; see e.g. \cite{_Griffi_Harri_}),
and $h:= h_0\cdot e^{-2\psi}$ the corresponding singular
metric. The curvature of $h$
is equal to $\6\bar\6\psi +\theta=\Xi$. 
We have shown that any nef bundle
admits a singular metric with positive current
as its curvature.

Let ${\cal I}$ denote the corresponding {\em multiplier ideal
sheaf}. It can be defined directly in terms of the function
$\psi$, but for our purposes it is more convenient to
define the tensor product $L\otimes {\cal I}$ directly
as a sheaf of all sections of $L$ which are locally
$L^2$-integrable in the singular metric $h$ defined 
above.

Assume now that $M$ is a compact K\"ahler
manifold with trivial
canonical bundle. By \ref{_Nadel_holo_coho_CY_Theorem_},
there is a natural surjection
\[ 
  H^0(\Omega^{n-i}\otimes L\otimes {\cal I}) \arrow H^i(L\otimes {\cal I}).
\]
To show that $M$ admits $L$-valued holomorphic
differential forms, it suffices to show that
$H^i(L)$ is non-vanishing, for some $i$. 

\hfill

\theorem\label{_coho_nonzero_then_effe_Theorem_}
Let $M$ be a simple hyperk\"ahler manifold, $\dim_\C M=2n$,
and  $L$ a nef bundle on $M$.
Consider a singular metric on $L$ with 
its curvature a positive current, 
and let ${\cal I}(L^m)$ be the sheaf of 
$L^2$-integrable holomorphic sections of $L^m$. 
\begin{description}
\item[(i)] Assume that for infinitely many  $m>0$,
$H^i({\cal I}(L^m))\neq 0$. Then
$L$ is $\Q$-effective.
\item[(ii)] Assume that all Lelong numbers of 
$L$ vanish. Then $L$ is $\Q$-effective.
\end{description}

\hfill

{\bf Proof:} By the multiplier ideal version
of Hard Lefschetz theorem \\ (\ref{_Nadel_holo_coho_CY_Theorem_}), 
$H^i({\cal I}(L^m))\neq 0$ implies existence of holomorphic
${\cal I}(L^m)$-valued differential forms on $M$.
However, ${\cal I}(L^m)$ is by construction 
a subsheaf of $L^m$, hence any ${\cal I}(L^m)$-valued differential form
can be considered as an $L^m$-valued differential form.
Now, \ref{_holo_L_valued_L_effe_Theorem_}
implies that $L^N$ is effective. This proves
\ref{_coho_nonzero_then_effe_Theorem_}  (i).
Then, \ref{_coho_nonzero_then_effe_Theorem_} (ii) follows,
because in this case  ${\cal I}(L^m)= L^m$, and
the usual calculation (see
\ref{_dim_V_k_cohomo_Corollary_}) implies that 
$\chi({\cal I}(L^m)=n+1$. \endproof

\hfill

{\bf Acknowledgements:}
I am grateful to E. Amerik, J.-P. Demailly,
G. Grantcharov, D. Kaledin, D. Matsushita, S. Nemirovsky, J. Sawon and
A. Todorov for interesting discussions. Many thanks to 
S\'ebastien  Boucksom for a valuable e-mail exchange.
Much gratitude to  F. Campana for interesting discussions
 and a reference to
 \cite{_Demailly_Peternell_Schneider:tan_nef_}.
Many thanks to the anonymous referees for numerous 
suggestions and remarks.
This paper was finished at IPMU (Tokyo University);
the author wishes to thank the IPMU staff for their
hospitality.

\hfill

\hfill

\small{

\noindent {\sc Misha Verbitsky\\
{\sc  Institute of Theoretical and
Experimental Physics \\
B. Cheremushkinskaya, 25, Moscow, 117259, Russia }\\
\tt verbit@mccme.ru }

}

\end{document}